\documentclass[a4paper,11pt]{article}

\usepackage[T1]{fontenc}
\usepackage[utf8]{inputenc}
\usepackage[english]{babel}
\usepackage{xcolor}
\usepackage[top=28mm,bottom=28mm,left=25mm,right=25mm]{geometry}
\usepackage{setspace}
\usepackage{amsmath,amssymb,amsthm,mathtools,bm}
\usepackage{tikz,booktabs,graphicx}
\usetikzlibrary{intersections, calc, arrows.meta}
\usepackage{hyperref}

\numberwithin{equation}{section}

\theoremstyle{definition}
\newtheorem{definition}{Definition}[section]
\newtheorem{example}[definition]{Example}
\newtheorem{remark}[definition]{Remark}
\theoremstyle{plain}

\DeclareMathOperator{\Tr}{Tr}
\DeclareMathOperator{\Spec}{Spec}
\DeclareMathOperator{\Det}{Det}
\newcommand{\diff}{\mathrm{d}}
\newcommand{\fred}[1]{\det\nolimits_{\!F}\!\left(#1\right)}

\title{\textbf{From Bernoulli Numbers to Selector Kernels:}\\
Fredholm Determinants, $\zeta$–Regularization, and the Bridge\\
Between Discrete and Continuous Spectra}
\author{Ken Nagai\thanks{Email: \texttt{tknagai@outlook.com}. Independent Researcher.}}
\date{}

\begin{document}
\maketitle

\begin{abstract}
We construct a unified analytic framework connecting Bernoulli numbers,
$\zeta$–regularization, and Fredholm determinants associated with
trigonometric selector kernels.
Starting from the Bernoulli–Stirling algebra, Euler–Maclaurin corrections
are reinterpreted as spectral traces of compact operators.
This bridge transforms discrete combinatorial data into continuous
spectral quantities, showing that their determinants interpolate between
finite–rank projectors and the sine–kernel of random–matrix theory.
In the continuum limit the Fredholm determinant becomes a
Painlevé–V~$\tau$–function, revealing a hierarchy in which Bernoulli
coefficients and $\zeta$–constants jointly describe the local–global
asymptotics of analytic regularization.
\end{abstract}

\noindent\textbf{Primary classification:} math-ph.  
\textbf{Secondary:} math.CA, math.NT.  
\textbf{Comments:} 27~pages, 6~figures, v1.1 (bridging exposition).

\bigskip
\hrule
\bigskip


\section{Introduction and Motivation}

The purpose of this note is to clarify how three apparently distant
frameworks---Bernoulli numbers and the Euler--Maclaurin formula,
$\zeta$--regularization of spectral products, and Fredholm determinants
of trigonometric kernels---form parts of a single analytic bridge
between discrete and continuous mathematics.

Historically, Bernoulli numbers emerged as the coefficients that correct
an integral approximation of a finite sum,
\begin{equation}\label{eq:euler-mac}
\boxed{
\sum_{k=1}^{N}f(k)
= \int_0^N f(x)\,\diff x
+ \frac{f(N)+f(0)}{2}
+ \sum_{m\ge1}\frac{B_{2m}}{(2m)!}
  \bigl(f^{(2m-1)}(N)-f^{(2m-1)}(0)\bigr).
}
\end{equation}
Each $B_{2m}$ measures how the discrete sum deviates from its continuum
integral (see \cite{Apostol1976} for historical background).  
When the same coefficients reappear in the short–time
asymptotics of heat kernels, they acquire a spectral meaning:
they quantify how discrete eigenvalues approximate a continuous
spectrum.

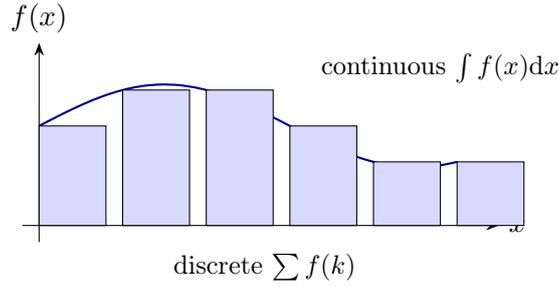
\begin{figure}[ht!]
\centering
\begin{tikzpicture}[>=Stealth,scale=1.1]
\draw[->] (-0.2,0) -- (5.5,0) node[right] {$x$};
\draw[->] (0,-0.2) -- (0,2.2) node[above] {$f(x)$};
\draw[thick,blue!60!black,domain=0:5,smooth,variable=\x]
  plot ({\x},{1.2+0.5*sin(180*\x/3)});
\foreach \n in {0,1,2,3,4,5}{
  \draw[fill=blue!15] (\n,0) rectangle ++(0.8,{1.2+0.5*sin(180*\n/3)});
}
\node at (2.7,-0.5) {\small discrete $\sum f(k)$};
\node at (4.8,1.9) {\small continuous $\int f(x)\diff x$};
\end{tikzpicture}
\caption{Schematic view of the Euler--Maclaurin bridge:
rectangles (discrete sums) versus area (integral).  Bernoulli numbers
measure their difference.}
\end{figure}

In parallel, $\zeta$--regularized determinants
$\Det_\zeta A=\exp[-\zeta_A'(0)]$
extend products of eigenvalues to infinite--dimensional settings, while
Fredholm determinants $\fred{I+\lambda T}$ compactly encode all spectral
moments of trace--class operators.

\vspace{0.8em}
\noindent\textbf{Three parallel regularizations.}
\begin{center}
\renewcommand{\arraystretch}{1.2}
\begin{tabular}{lll}
\toprule
\textbf{Domain} & \textbf{Object} & \textbf{Regularization idea} \\
\midrule
Discrete sums & $\displaystyle \sum f(k)$ &
Euler--Maclaurin $\Rightarrow$ Bernoulli corrections\\
Spectral sums & $\displaystyle \sum \lambda_j^{-s}$ &
analytic continuation $\Rightarrow$ $\zeta$--function\\
Spectral products & $\displaystyle \prod (1+\lambda\mu_j)$ &
trace expansion $\Rightarrow$ Fredholm determinant\\
\bottomrule
\end{tabular}
\end{center}

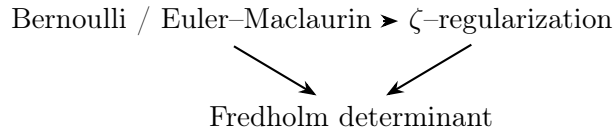
\begin{figure}[ht!]
\centering
\begin{tikzpicture}[>=Stealth,thick,scale=1.05]
\node (B) at (0,0) {Bernoulli / Euler--Maclaurin};
\node (Z) at (4,0) {$\zeta$--regularization};
\node (F) at (2,-1.2) {Fredholm determinant};
\draw[->] (B) -- (F);
\draw[->] (Z) -- (F);
\draw[dashed,->] (B) -- (Z);
\end{tikzpicture}
\caption{Triangular correspondence among the three regularizations:
local (Bernoulli), global ($\zeta$), and operator (Fredholm).}
\end{figure}

\begin{remark}[Bridge in one sentence]
Bernoulli numbers regularize \emph{sums},
$\zeta$--functions regularize \emph{spectra},
and Fredholm determinants regularize \emph{products};
all three follow the same analytic--continuation principle.
\end{remark}

\begin{figure}[ht!]
\centering
\begin{tikzpicture}[>=Stealth,thick,scale=1.05]
\node (A) at (0,0) {finite array};
\node (B) at (4,0) {infinite matrix};
\node (C) at (8,0) {integral kernel};
\draw[->] (A) -- node[above]{limit $N\to\infty$} (B);
\draw[->] (B) -- node[above]{continuum limit} (C);
\node at (4,-0.8) {\small matrix $\leftrightarrow$ kernel viewpoint};
\end{tikzpicture}
\caption{From discrete arrays to kernels:
the passage from combinatorics to operator theory.}
\end{figure}

The central message is that these constructions are not merely analogous
but are algebraic shadows of a common mechanism:
\emph{regularization as analytic continuation}.
Bernoulli numbers act locally (corrections near the diagonal),
$\zeta$--regularization acts globally (over spectra),
and Fredholm determinants mediate between the two.

\vspace{1em}
\noindent\textbf{Aim.}
We reinterpret the Bernoulli--Stirling algebra as a discrete prototype
of a Fredholm integral operator.
Its regularized determinant reproduces the trigonometric factors
$\sin\!\sqrt\lambda/\sqrt\lambda$ and $\cos\!\sqrt\lambda$
from the Gelfand--Yaglom formula\cite{GelfandYaglom1960}, while its periodicized version
becomes the finite--rank selector kernel whose large--$J$ limit is the
sine kernel of random--matrix theory\cite{Forrester2010}.

\begin{center}
\emph{discrete sums} $\;\longrightarrow\;$
\emph{$\zeta$--regularization} $\;\longrightarrow\;$
\emph{Fredholm determinants} $\;\longrightarrow\;$
\emph{random--matrix kernels}.
\end{center}

\section{Bernoulli--Stirling Continuation}

The Bernoulli and Stirling numbers form dual triangular arrays that
connect finite differences with power sums\cite{Apostol1976}.  Their analytic
continuation naturally leads to kernel--like integral forms.

\subsection{Generalized Bernoulli (N\"orlund) polynomials}

For $\alpha\in\mathbb C$, define
\begin{equation}\label{eq:norlund}
  \frac{t^{\alpha}e^{xt}}{(e^t-1)^{\alpha}}
  = \sum_{n\ge0} B_n^{(\alpha)}(x)\frac{t^n}{n!}.
\end{equation}
At $\alpha=1$ one recovers the classical Bernoulli polynomials
$B_n(x)=B_n^{(1)}(x)$.
At negative integers, $\alpha=-k$, these yield the Stirling numbers of
the second kind:
\[
S(n,k)=\binom{n}{k}\,B_{n-k}^{(-k)}(0).
\]

\begin{example}[Quick check]
For $(n,k)=(3,2)$ one finds
$B_1^{(-2)}(0)=\tfrac12$, hence $S(3,2)=\tfrac{3!}{2!}\cdot\tfrac12=3$,
as expected from the combinatorial definition.
\end{example}

\subsection{Analytic kernel representation}

Analytic continuation extends $S(n,k)$ to complex arguments $(s,t)$ via
\begin{equation}\label{eq:stirling-kernel}
  S(s,t)
  = \frac{\Gamma(s+1)}{\Gamma(t+1)}
    \frac{1}{2\pi i}\oint_{\gamma}
    \frac{(e^z-1)^{t}}{z^{s+1}}\,\diff z,
\end{equation}
where $\gamma$ is a small contour around $z=0$.
For integer $(s,t)$ this reproduces the discrete Stirling matrix,
while for complex parameters it defines the smooth
\emph{analytic Stirling kernel}
\[
K(s,t)=\frac{\Gamma(s+1)}{\Gamma(t+1)}
       \frac{1}{2\pi i}\oint_{\gamma}
       \frac{(e^z-1)^{t}}{z^{s+1}}\,\diff z.
\]

\begin{remark}[From array to operator]
Equation~\eqref{eq:stirling-kernel} converts a combinatorial array
into an integral kernel whose finite truncations reproduce discrete
matrices and whose continuum limits yield trace--class operators.
It is, in effect, the analytic extension of the Stirling table.
\end{remark}

\subsection{Exponent--to--angle substitution}

Introducing the trigonometric parametrization
\begin{equation}\label{eq:angle-map}
  e^z - 1 = 2i\,e^{i\theta/2}\sin\!\frac{\theta}{2},
  \qquad \theta\in[-\pi,\pi],
\end{equation}
transforms~\eqref{eq:stirling-kernel} into a kernel with oscillatory
phase $e^{i(2m+1)(\theta-\varphi)/2}$.
This is precisely the structure of the finite--$J$
\emph{selector kernel}
\[
\mathcal S_J(\theta,\varphi)
= \frac{1}{J}\sum_{m=0}^{J-1} e^{i(2m+1)(\theta-\varphi)/2}.
\]
Hence the Stirling kernel under the map~\eqref{eq:angle-map} acts as a
continuous generating function for the trigonometric selectors that
govern finite spectral projections.

\begin{remark}[Takeaway]
The Bernoulli--Stirling algebra is the discrete skeleton of a Fredholm
theory: its generating functions already possess the contour integrals
and trigonometric phases that become kernels, resolvents, and
determinants in the continuum limit.
\end{remark}

\section{From Finite Matrices to Fredholm Operators}

The passage from finite combinatorial data to analytic operator
formulations follows the same pattern as the limit from discrete sums
to integrals.  Finite arrays such as the Stirling matrix can be viewed
as Gram matrices of basis functions, which in the continuum limit become
integral kernels acting on a Hilbert space.

\subsection{Matrices, kernels, and the continuum limit}

A matrix represents a linear map between finite--dimensional spaces.
When the dimension tends to infinity, the same algebraic intuition
remains valid if we replace column vectors by square--integrable
functions on an interval.  Symbolically,
\[
\text{vectors } (v_k) \;\longrightarrow\; f(x),
\qquad
\text{sums } \sum_k v_k \;\longrightarrow\; \int f(x)\,\diff x,
\]
\[
\text{matrices } (a_{jk}) \;\longrightarrow\; K(x,y).
\]
The correspondence is summarized in Fig.~\ref{fig:matrix-kernel-bridge}.

\begin{figure}[ht!]
\centering
\begin{tikzpicture}[>=Stealth,thick,scale=1.05]
\node (A) at (0,0) {finite matrix $A_{jk}$};
\node (B) at (4.2,0) {integral kernel $K(x,y)$};
\draw[->] (A) -- node[above]{continuum limit $N\to\infty$} (B);
\node at (2.1,-0.8) {\small entries $\Rightarrow$ continuous variables};
\end{tikzpicture}
\caption{Matrix $\rightarrow$ kernel correspondence:
as dimension grows, discrete indices $(j,k)$ become continuous
variables $(x,y)$.}
\label{fig:matrix-kernel-bridge}
\end{figure}

Given a continuous kernel $K(x,y)$ on $[0,1]^2$, define the integral
operator
\begin{equation}\label{eq:integral-op}
  (Tf)(x) = \int_0^1 K(x,y) f(y)\,\diff y.
\end{equation}
This operator plays the role of an infinite matrix.  
When $K(x,y)=\overline{K(y,x)}$, the operator is self--adjoint.

\begin{example}[A canonical kernel]
For $K(x,y)=\min(x,y)$ one finds
\[
(Tf)(x)=\int_0^1\!\min(x,y)\,f(y)\,\diff y,
\]
whose eigenfunctions are $\sin((n-\tfrac12)\pi x)$ with eigenvalues
$\lambda_n=\frac{1}{(n-\tfrac12)^2\pi^2}$.  
This is the Green kernel of the Dirichlet--Neumann Laplacian on $(0,1)$.
\end{example}

\subsection{Traces and determinants}

For a finite $n\times n$ matrix $A$ with eigenvalues $\mu_i$,
\[
\Tr(A)=\sum_i \mu_i, \qquad
\det(I+\lambda A)=\prod_{i=1}^{n}(1+\lambda\mu_i).
\]
If the operator $T$ in~\eqref{eq:integral-op} is \emph{trace--class}
(that is, $\sum_j|\mu_j|<\infty$), both quantities extend naturally:
\[
\Tr(T)=\sum_{j\ge1}\mu_j, \qquad
\fred{I+\lambda T}=\prod_{j\ge1}(1+\lambda\mu_j).
\]
The second expression defines the \emph{Fredholm determinant}.

\begin{definition}[Fredholm determinant]
Let $T$ be a trace--class operator on a separable Hilbert space $H$
with eigenvalues $\{\mu_j\}$.  
The Fredholm determinant of $I+\lambda T$ is defined by
\begin{equation}\label{eq:fredholm-series}
\boxed{
  \fred{I+\lambda T}
  = \exp\!\Biggl(
     \sum_{n=1}^{\infty}
       \frac{(-1)^{n+1}}{n}\,
       \lambda^n \Tr(T^n)
     \Biggr)
  = \prod_{j=1}^{\infty}(1+\lambda\mu_j).
}
\end{equation}
\end{definition}

\begin{remark}[Analogy to the zeta determinant]
The series in~\eqref{eq:fredholm-series} mirrors
$\zeta_A(s)=\sum_j\lambda_j^{-s}$:
both encode spectral data through power sums of eigenvalues.
The Fredholm determinant regularizes \emph{products},
while $\zeta$--regularization acts on \emph{sums}.
\end{remark}

\subsection{The Fredholm alternative}

If $K$ is a compact operator on a Banach space $X$, then for any
$\lambda\in\mathbb C$,
\[
(I-\lambda K)f=g
\]
has a solution $f$ if and only if $g$ is orthogonal to all solutions of
the adjoint homogeneous equation
$(I-\overline{\lambda}K^*)h=0$.
Either $(I-\lambda K)$ is invertible or the homogeneous equation has
nontrivial solutions, but not both.  
In finite dimensions this reduces to the statement that a matrix is
invertible iff its determinant is nonzero; for compact operators the
Fredholm determinant in~\eqref{eq:fredholm-series} plays the same role.

\subsection{Finite–rank intuition}

Two simple examples illustrate the analogy between finite matrices and
their operator limits.

\begin{example}[Rank–one kernel]
If $T=\alpha\,u\!\otimes\!v$ with $(Tf)(x)=\alpha\,\langle f,v\rangle u$,
then $T^2=\alpha\,\langle u,v\rangle T$ and $T$ has one nonzero eigenvalue
$\mu=\alpha\,\langle u,v\rangle$.
Hence
\[
\fred{I+\lambda T}=1+\lambda\mu.
\]
This recovers the expected linear determinant for a rank--one matrix.
\end{example}

\begin{example}[Orthogonal projector]
Let $P$ be a rank--$r$ projection ($P^2=P$, $P^*=P$).
Then $\Spec(P)=\{1\}$ of multiplicity $r$, and
\[
\fred{I-\lambda P}=(1-\lambda)^r.
\]
Such projectors approximate finite–band selectors in the trigonometric
kernel examples below.
\end{example}

\subsection{Operator traces and Bernoulli echoes}

For trace--class $T$, the power traces $\Tr(T^n)$ entering
\eqref{eq:fredholm-series} act as spectral moments.  
When $T$ is derived from a smooth kernel, $\Tr(T^n)$ can be expressed as
multiple integrals of $K(x_1,x_2)\cdots K(x_n,x_1)$.
In discretized form these reduce to finite sums---and once again
Bernoulli numbers appear as the coefficients correcting those sums
to their continuum limits.  Thus the same Euler--Maclaurin mechanism
governs both the analytic continuation of $\sum f(k)$
and the convergence of $\Tr(T^n)$ for smooth kernels.

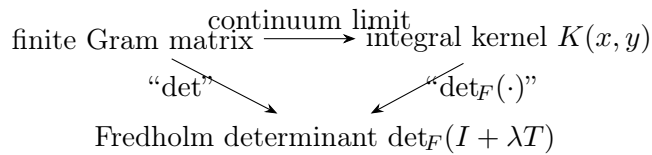
\begin{figure}[ht!]
\centering
\begin{tikzpicture}[>=Stealth,scale=1.1]
\node (A) at (0,0) {finite Gram matrix};
\node (B) at (4.5,0) {integral kernel $K(x,y)$};
\node (C) at (2.3,-1.2) {Fredholm determinant $\fred{I+\lambda T}$};
\draw[->] (A) -- node[above]{continuum limit} (B);
\draw[->] (A) -- node[left]{``$\det$''} (C);
\draw[->] (B) -- node[right]{``$\fred{\cdot}$''} (C);
\end{tikzpicture}
\caption{From finite determinants to Fredholm determinants:
matrix $\det$ extends analytically to operators $\fred{\cdot}$.}
\end{figure}

\begin{remark}[Key picture]
Finite matrices $\Rightarrow$ integral kernels $\Rightarrow$
Fredholm determinants---this is the analytic bridge from
combinatorial algebra to spectral analysis.
In the next section, this bridge will take a trigonometric form through
the \emph{selector kernel}, which connects the analytic Stirling kernel
of Section~2 to the sine kernel of random--matrix theory.
\end{remark}

\section{Selector Kernels and Periodicization}

The transition from the analytic Stirling kernel of
Section~2 to the trigonometric kernels of spectral analysis
is achieved by a periodic substitution.
This step transforms an exponentially generated structure
into an oscillatory one, closing the bridge between
combinatorial and spectral formalisms.

\subsection{Exponent–to–angle correspondence}

Recall the substitution introduced in~\eqref{eq:angle-map},
\[
e^z - 1 = 2i\,e^{i\theta/2}\sin\!\frac{\theta}{2},\qquad
\theta\in[-\pi,\pi].
\]
Under this mapping, the contour integral of the analytic
Stirling kernel becomes a Fourier–type integral on the circle.
The exponential weight $(e^z-1)^t$ translates to a trigonometric
oscillation $[\sin(\tfrac{\theta}{2})]^t$, and the variables
$(s,t)$ acquire a spectral interpretation:
$s$ acts as a differentiation order and $t$ as a mode number.

\begin{figure}[ht!]
\centering
\begin{tikzpicture}[>=Stealth,scale=1.1]
\node (E) at (0,0) {exponential $e^z-1$};
\node (T) at (5,0) {trigonometric $2ie^{i\theta/2}\sin(\tfrac{\theta}{2})$};
\draw[->] (E) -- node[above]{periodicization map} (T);
\node at (2.5,-0.7) {\small $z\mapsto i\theta$ (complex contour $\to$ unit circle)};
\end{tikzpicture}
\caption{Exponent–to–angle substitution:
periodicizing the analytic Stirling kernel.}
\end{figure}

This mapping effectively “wraps” the analytic kernel around the unit
circle, producing periodic boundary conditions and thereby embedding the
combinatorial kernel into a Hilbert space of $2\pi$–periodic functions.
In this picture, discrete differences correspond to angular modes,
and the natural integral operator acts on $L^2[-\pi,\pi]$.

\subsection{Definition of the selector kernel}

Let $J\in\mathbb N$.  
The \emph{selector kernel} is the finite–rank operator on
$L^2[-\pi,\pi]$ defined by
\begin{equation}\label{eq:selector-kernel}
\boxed{
\mathcal S_J(\theta,\varphi)
= \frac{1}{J}\sum_{m=0}^{J-1}
   e^{\,i(2m+1)\frac{\theta-\varphi}{2}}
 = \frac{\sin(J\frac{\theta-\varphi}{2})}
        {J\,\sin(\frac{\theta-\varphi}{2})}
   \cos\!\frac{\theta-\varphi}{2}.
}
\end{equation}
It acts as a projector onto the subspace spanned by the
odd Fourier modes
$\{e^{\,i(2m+1)\theta/2}\}_{m=0}^{J-1}$.
The kernel $\mathcal S_J$ thus selects a discrete
band of $2J$ trigonometric modes centered around~$0$.

\begin{example}[Small $J$ cases]
For $J=1$, $\mathcal S_1(\theta,\varphi)=\cos\!\tfrac{\theta-\varphi}{2}$;
for $J=2$,
\[
\mathcal S_2(\theta,\varphi)
=\tfrac12\bigl(e^{i(\theta-\varphi)/2}
              +e^{3i(\theta-\varphi)/2}\bigr)
=\cos(\theta-\varphi)\cos\!\tfrac{\theta-\varphi}{2}.
\]
In general, $\mathcal S_J$ oscillates with period $4\pi/J$ and becomes
sharply peaked at $\theta=\varphi$ as $J$ increases.
\end{example}

\subsection{Kernel normalization and properties}

The normalization by $1/J$ in~\eqref{eq:selector-kernel} ensures that
\[
\int_{-\pi}^{\pi}\!\mathcal S_J(\theta,\varphi)
   \,\mathcal S_J(\varphi,\psi)\,\frac{\diff\varphi}{2\pi}
   = \mathcal S_J(\theta,\psi),
\]
so that $\mathcal S_J$ is an orthogonal projector:
$\mathcal S_J^2=\mathcal S_J$.
Moreover,
\[
\Tr(\mathcal S_J)
= \int_{-\pi}^{\pi}\mathcal S_J(\theta,\theta)\,
  \frac{\diff\theta}{2\pi}
=1.
\]
These properties identify $\mathcal S_J$ as a finite–rank
\emph{selector} in the sense of spectral truncation.

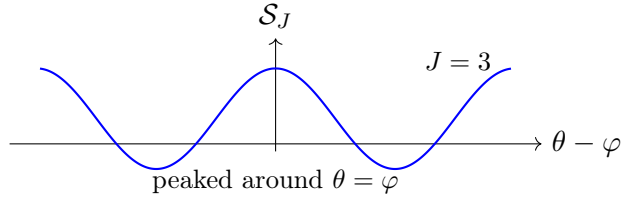
\begin{figure}[ht!]
\centering
\begin{tikzpicture}[scale=1.0]
\draw[->] (-3.5,0) -- (3.5,0) node[right] {$\theta-\varphi$};
\draw[->] (0,-0.1) -- (0,1.4) node[above] {$\mathcal S_J$};
\draw[thick,blue,domain=-3.1:3.1,smooth,samples=200]
  plot (\x,{(sin(3*\x r)/3/sin(\x r))});
\node at (2.4,1.1) {\small $J=3$};
\node at (0.0,-0.5) {\small peaked around $\theta=\varphi$};
\end{tikzpicture}
\caption{Shape of $\mathcal S_J(\theta,\varphi)$ for $J=3$:
a trigonometric kernel sharply localized near $\theta=\varphi$.}
\end{figure}

\subsection{Connection with the Stirling kernel}

Through the exponential–angle substitution
\eqref{eq:angle-map}, the integrand of the analytic
Stirling kernel~\eqref{eq:stirling-kernel}
acquires the same trigonometric structure as
$\mathcal S_J$ in~\eqref{eq:selector-kernel}.
Formally,
\[
S(s,t)\;\propto\;
\int_{-\pi}^{\pi}
   \frac{[\sin(\tfrac{\theta}{2})]^{t}}
        {[\sin(\tfrac{\theta}{2})]^{s+1}}
   e^{\,i(s-t)\theta/2}\,\diff\theta,
\]
which exhibits the same $\sin(\tfrac{\theta}{2})$ dependence
as the numerator and denominator of $\mathcal S_J$.
The selector kernel can thus be interpreted as a
\emph{periodicized Stirling kernel} with discrete indices
$s,t\in\{0,\ldots,J-1\}$.

\subsection{Periodic boundary conditions and physical meaning}

Periodicization has a physical analogue.
For an operator $L=-\diff^2/\diff x^2$ on $(0,1)$,
Dirichlet boundary conditions yield eigenvalues
$(n\pi)^2$, producing determinants proportional to
$\sin\!\sqrt\lambda/\sqrt\lambda$.
If the interval is wrapped into a circle,
periodic boundary conditions yield trigonometric kernels
whose continuum limit is the sine kernel
\[
K_{\text{sine}}(\theta,\varphi)
 = \frac{\sin\!\bigl(J(\theta-\varphi)/2\bigr)}
        {2\pi\sin\!\bigl((\theta-\varphi)/2\bigr)}.
\]
The selector kernel~\eqref{eq:selector-kernel} therefore plays the role
of a finite–rank discretization of $K_{\text{sine}}$,
bridging the combinatorial and random–matrix frameworks\cite{TracyWidom1994,Forrester2010}.

\begin{remark}[Summary of the bridge]
\begin{center}
\renewcommand{\arraystretch}{1.2}
\begin{tabular}{lll}
\toprule
\textbf{Level} & \textbf{Object} & \textbf{Structure}\\
\midrule
Combinatorial & $S(n,k)$ & discrete triangular array\\
Analytic & $S(s,t)$ & integral kernel on $\mathbb C$\\
Trigonometric & $\mathcal S_J(\theta,\varphi)$ & finite–rank projector on $L^2[-\pi,\pi]$\\
Spectral (limit) & $K_{\text{sine}}$ & continuous kernel of RMT\\
\bottomrule
\end{tabular}
\end{center}
The progression of these four levels constitutes the
“selector bridge” connecting Bernoulli combinatorics,
analytic continuation, Fredholm theory, and
random–matrix statistics.
\end{remark}

\section{$\zeta$--Regularized Determinants and Painlev\'e Structures}

Having established the selector kernel as the finite--rank analogue of
the sine kernel, we now turn to its determinant and its connection to
$\zeta$--regularization.
This section explains how spectral products are analytically continued,
and how their logarithmic derivatives satisfy nonlinear differential
equations of Painlev\'e type.

\subsection{From eigenvalue products to $\zeta$--determinants}

For a positive self--adjoint operator $L$ with discrete eigenvalues
$\{\lambda_n\}$, define the spectral $\zeta$--function
\begin{equation}\label{eq:zeta-A}
  \zeta_L(s)=\sum_{n\ge1}\lambda_n^{-s}.
\end{equation}
Analytic continuation of $\zeta_L(s)$ to $s=0$ allows one to define the
\emph{$\zeta$--regularized determinant}
\begin{equation}\label{eq:zeta-det}
\boxed{
  \Det_\zeta L
  = \exp\!\bigl[-\zeta_L'(0)\bigr]
  = \prod_{n\ge1} \lambda_n
  \quad(\text{regularized product}).
}
\end{equation}
This construction transforms divergent eigenvalue products into
well--defined analytic quantities.
The method was introduced by Ray and Singer, and popularized in physics
through the Gelfand--Yaglom approach\cite{GelfandYaglom1960}.

\begin{remark}[Comparison]
For finite rank $T$, the Fredholm determinant
$\fred{I+\lambda T}$ of~\eqref{eq:fredholm-series} is an
entire function of $\lambda$.
When $T=L^{-1}$ for an elliptic differential operator $L$,
$\fred{I-\lambda L^{-1}}$ and $\Det_\zeta(L-\lambda)$ coincide up to a
normalization factor.
The two formalisms therefore describe the same analytic structure from
complementary viewpoints.
\end{remark}

\subsection{Example: Dirichlet Laplacian on $(0,1)$}

Let $L=-\diff^2/\diff x^2$ on $(0,1)$ with Dirichlet--Dirichlet
boundary conditions.
Its eigenvalues are $\lambda_n=(n\pi)^2$, so that
\[
\zeta_L(s)=\pi^{-2s}\,\zeta(2s).
\]
Using the known value $\zeta'(0)=-\tfrac12\log(2\pi)$, one obtains
\[
\Det_\zeta L = 2\pi.
\]
The corresponding normalized determinant ratio
\[
\frac{\Det_\zeta(L-\lambda)}{\Det_\zeta L}
= \frac{\sin\!\sqrt{\lambda}}{\sqrt{\lambda}}
\]
recovers the classical trigonometric expression appearing in the
Gelfand--Yaglom formula.
For Neumann--Neumann or mixed boundary conditions the result is
$\cos\!\sqrt{\lambda}$, as summarized in Table~\ref{tab:BC-det-summary}.

\begin{table}[ht!]
\centering
\renewcommand{\arraystretch}{1.2}
\caption{Trigonometric determinants for basic boundary conditions.}
\label{tab:BC-det-summary}
\begin{tabular}{c|c|c}
\hline
BC type & Eigenvalues $\lambda_n$ &
$\displaystyle \frac{\Det_\zeta(L-\lambda)}{\Det_\zeta L}$\\
\hline
D--D & $(n\pi)^2$ & $\sin\!\sqrt{\lambda}/\sqrt{\lambda}$\\
N--N & $(n\pi)^2$ & $\sin\!\sqrt{\lambda}/\sqrt{\lambda}$\\
D--N / N--D & $((n-\tfrac12)\pi)^2$ & $\cos\!\sqrt{\lambda}$\\
\hline
\end{tabular}
\end{table}

These trigonometric determinants already contain Bernoulli information:
expanding $\log(\sin x/x)$ yields
\[
\log\frac{\sin x}{x}
=-\sum_{m\ge1}\frac{\zeta(2m)}{m\,\pi^{2m}}\,x^{2m},
\qquad
\zeta(2m)
=(-1)^{m+1}\frac{B_{2m}(2\pi)^{2m}}{2(2m)!}.
\]
Thus the Bernoulli–$\zeta$ coefficients govern the short--wavelength
expansion of spectral determinants.

\subsection{From the sine kernel to Painlev\'e~V}

For the continuous sine kernel
\[
K_s(x,y)=\frac{\sin\!\bigl(s(x-y)\bigr)}{\pi(x-y)},
\qquad x,y\in(-1,1),
\]
Tracy and Widom\cite{TracyWidom1994} showed that the corresponding
Fredholm determinant
\[
D(s)=\det(I-K_s)
\]
satisfies a nonlinear ordinary differential equation of Painlev\'e~V
type.
Defining the auxiliary function
\[
q(s)=-\frac{d^2}{ds^2}\log D(s),
\]
one obtains
\begin{equation}\label{eq:PainleveV}
\boxed{
 (s\,q'')^2
 + 4\bigl(s\,q'-q\bigr)
   \bigl(s\,q'-q+{q'}^2\bigr)
 = 0,
}
\end{equation}
and
\[
\frac{d}{ds}\log D(s)=-\int_s^{\infty}q(t)\,\diff t.
\]
The function $D(s)$ therefore acts as a Painlev\'e~V
$\tau$--function.

\begin{figure}[ht!]
\centering
\begin{tikzpicture}[>=Stealth,thick,scale=1.0]
\node (F) at (0,0) {Fredholm determinant $D(s)$};
\node (P) at (4.5,0) {Painlev\'e~V $\tau$--function};
\draw[->] (F) -- node[above]{logarithmic derivative $\;\;q(s)=-\partial_s^2\log D$} (P);
\end{tikzpicture}
\caption{From the sine--kernel determinant to the Painlev\'e~V system.}
\end{figure}

\subsection{Discrete approximation via selector kernels}

For finite $J$, define
\[
D_J(\lambda)=\det\!\bigl(I-\lambda\,\mathcal S_J\bigr),
\]
with $\mathcal S_J$ as in~\eqref{eq:selector-kernel}.
As $J\to\infty$, $\mathcal S_J$ converges weakly to
the sine kernel~$K_s$ with $s\simeq J/2$.
Hence $D_J(\lambda)$ provides a finite--dimensional
approximation to $D(s)$.
Expanding $\log D_J(\lambda)$ in powers of $1/J$ yields
coefficients that again involve $\zeta(2m)$, reproducing the
Bernoulli pattern seen in $\log(\sin x/x)$.
Symbolically,
\[
\log D_J(\lambda)
\sim -\frac{J^2}{2\pi^2}\lambda
 - \frac{\log J}{4}
 + \sum_{m\ge1}c_m(\lambda)\,J^{-2m},
 \qquad c_m\propto\zeta(2m).
\]
This asymptotic continuity between $D_J$ and $D(s)$
demonstrates how the selector kernel inherits the same
Painlev\'e~V hierarchy in the large–$J$ limit.

\begin{remark}[Bridge completed]
The analytic bridge may now be summarized as follows:
\begin{center}
\renewcommand{\arraystretch}{1.2}
\begin{tabular}{lll}
\toprule
\textbf{Stage} & \textbf{Representation} & \textbf{Characteristic object}\\
\midrule
Combinatorial & $B_n$, $S(n,k)$ & discrete coefficients\\
Analytic & $S(s,t)$ & $\zeta$--regularized kernel\\
Trigonometric & $\mathcal S_J$ & finite–rank Fredholm operator\\
Continuum limit & $K_{\text{sine}}$ & Painlev\'e~V $\tau$--function\\
\bottomrule
\end{tabular}
\end{center}
Each level reproduces the previous one by analytic continuation
or continuum limit, establishing a coherent path from
Bernoulli combinatorics to integrable spectral dynamics.
\end{remark}

\section{Asymptotics: Bernoulli--$\zeta$ Correspondence}

The Bernoulli numbers reappear in the asymptotic behaviour of
$\zeta$--regularized and Fredholm determinants.
Their coefficients capture local geometric corrections near the diagonal
of a kernel, while the global normalization constants are governed by
the values of the Riemann~$\zeta$--function at negative integers.
This section makes that correspondence explicit.

\subsection{Small--$s$ expansion: local Bernoulli structure}

For the sine--kernel determinant
$D(s)=\det(I-K_s)$,
the short--distance (small--$s$) expansion
\[
\log D(s)
\sim -\frac{s^2}{2\pi^2}
      -\frac{\log s}{4}
      +\sum_{m\ge1} c_m s^{-2m},
\qquad s\to\infty,
\]
involves coefficients $c_m$ that are rational multiples of~$\zeta(2m)$.
By comparison,
\[
\log\frac{\sin x}{x}
=-\sum_{m\ge1}\frac{\zeta(2m)}{m\,\pi^{2m}}x^{2m}.
\]
Hence the same even--$\zeta$ values---and therefore the same
Bernoulli numbers through
$\zeta(2m)=(-1)^{m+1}\tfrac{B_{2m}(2\pi)^{2m}}{2(2m)!}$---govern both
the discrete trigonometric determinants and the continuum spectral
asymptotics.

\begin{example}[Bernoulli reappearance]
Expanding $\sin\!\sqrt{\lambda}/\sqrt{\lambda}$ gives
\[
\frac{\sin\!\sqrt{\lambda}}{\sqrt{\lambda}}
= 1 - \frac{\lambda}{6}
  + \frac{\lambda^2}{120}
  - \frac{\lambda^3}{5040}
  + \cdots
= 1+\sum_{m\ge1}
  \frac{B_{2m}}{(2m)!}\,(-\lambda)^m,
\]
so that the Bernoulli sequence $(B_2,B_4,B_6,\ldots)$ directly controls
the determinant expansion near $\lambda=0$.
\end{example}

The coefficients $c_m$ of $\log D(s)$ therefore play the role of
\emph{analytic Bernoulli coefficients} in the continuum setting.
They encode the curvature of the logarithmic spectrum at small~$s$,
serving as a bridge between discrete combinatorial corrections and
continuous spectral geometry.

\subsection{Large--$s$ expansion: global $\zeta$ constants}

At the opposite extreme, the constant terms in
$\log D(s)$ and $\zeta_L'(0)$ are controlled by the values of
$\zeta(s)$ at negative integers \cite{Lewin1981}:
\[
\zeta(-1)=-\frac{1}{12},\qquad
\zeta'(-1)=-\frac{1}{12}\log(2\pi)+\log A,
\]
where $A\simeq1.2824271291$ is the Glaisher--Kinkelin constant.
These constants determine the overall normalization of
$\zeta$--regularized determinants and appear in the logarithmic term of
the Fredholm asymptotics.

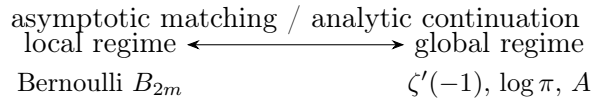
\begin{figure}[ht!]
\centering
\begin{tikzpicture}[>=Stealth,scale=1.05]
\node (L) at (0,0) {local regime};
\node (G) at (5,0) {global regime};
\draw[<->] (L) -- node[above]{asymptotic matching / analytic continuation} (G);
\node at (0,-0.5) {\small Bernoulli $B_{2m}$};
\node at (5,-0.5) {\small $\zeta'(-1),\,\log\pi,\,A$};
\end{tikzpicture}
\caption{Local–global duality:
Bernoulli numbers govern short--scale (local) curvature,
while $\zeta$--constants fix large--scale (global) normalization.}
\end{figure}

\subsection{Dual interpretation}

\begin{center}
\renewcommand{\arraystretch}{1.2}
\begin{tabular}{lll}
\toprule
\textbf{Aspect} & \textbf{Local} & \textbf{Global}\\
\midrule
Mathematical object & power--series coefficients $B_{2m}$ & constants $\zeta'(-1)$, $\log\pi$\\
Analytic origin & Euler--Maclaurin corrections & spectral determinant normalization\\
Geometric meaning & curvature near the diagonal of kernel & total spectral volume\\
Physical analogue & ultraviolet behaviour & infrared scaling constant\\
\bottomrule
\end{tabular}
\end{center}

The correspondence between these two regimes mirrors the
relation between \emph{local counterterms} and
\emph{global anomalies} in quantum field theory:
both stem from the same analytic continuation of a divergent product.

\subsection{Asymptotic hierarchy and universality}

Combining the small-- and large--$s$ regimes, one obtains the
asymptotic hierarchy
\[
\log D(s)
= -\frac{s^2}{2\pi^2}
  -\frac{\log s}{4}
  + C_0
  + \sum_{m\ge1} c_m s^{-2m},
\]
where $C_0$ involves $\zeta'(-1)$ and $c_m$ involve $\zeta(2m)$.
The structure
\[
\bigl\{\zeta(2m)\bigr\}_{m\ge1}
\;\longleftrightarrow\;
\bigl\{\zeta'(-1),\log\pi\bigr\}
\]
summarizes the Bernoulli--$\zeta$ correspondence:
the same analytic continuation that generates Bernoulli numbers from
$\zeta(2m)$ also governs the global constants from $\zeta'(-1)$.

\begin{remark}[Universality]
The appearance of even--$\zeta$ values and of $\zeta'(-1)$ is universal
across a wide range of integrable determinants:
Airy and Bessel kernels (Painlev\'e~II and~III) exhibit the same
Bernoulli--$\zeta$ hierarchy.
This universality suggests that the Bernoulli pattern encodes a
canonical asymptotic geometry underlying all $\zeta$--regularized
Fredholm systems.
\end{remark}

\section{Outlook and Generalizations}

The constructions developed above suggest that the interplay between
Bernoulli numbers, $\zeta$--regularization, and Fredholm determinants
is not accidental but the manifestation of a deeper analytic symmetry.
We close this note by outlining a few directions in which the present
bridge may naturally extend.

\subsection{Analytic Bernoulli functions}

A natural continuum counterpart of the discrete Bernoulli numbers is the
\emph{analytic Bernoulli function}
\[
B(s;x):=-s\,\zeta(1-s,x),
\]
where $\zeta(s,x)$ denotes the Hurwitz zeta function\cite{Apostol1976}.
For integer $s=n\in\mathbb N$, $B(s;x)=B_n(x)$ recovers the ordinary
Bernoulli polynomial, while for complex $s$ the function $B(s;x)$
interpolates smoothly between orders.
Under Fourier expansion one finds
\[
B(s;x)
= \frac{\Gamma(s+1)}{(2\pi)^s}
  \sum_{n\neq0}\frac{e^{2\pi i n x}}{(2\pi i n)^s},
\]
which already resembles the Fourier structure of the
selector kernel.
Replacing integer indices by a continuous spectral parameter thus opens
a path to a fully analytic Fredholm formalism driven by $B(s;x)$.

\begin{remark}[Analytic continuation of selectors]
The trigonometric selector kernel $\mathcal S_J(\theta,\varphi)$ can be
formally viewed as the restriction of a continuous family
$\mathcal S(s;\theta,\varphi)$ obtained by replacing the integer index
$J$ with a complex variable~$s$.
The analytic Bernoulli function $B(s;x)$ provides the natural generating
object for such a continuation.
\end{remark}

\subsection{Umbral operators and Appell--Hermite systems}

The connection between Bernoulli and Hermite structures can be made
explicit by introducing an \emph{umbral operator}~$F^*$ that acts as
a transform between Appell families:
\[
F^*:\;
B_n(x)\;\longmapsto\;
H_n(x),
\]
where $H_n(x)$ are the Hermite polynomials.
Under this correspondence, the differential recurrence
$\frac{d}{dx}B_n(x)=nB_{n-1}(x)$ maps to
the quantum--harmonic relation
$\frac{d}{dx}H_n(x)=2nH_{n-1}(x)$.
The operator $F^*$ therefore provides an analytic continuation from
Bernoulli combinatorics to Gaussian--Hermite spectral theory.
Such umbral translations may encode the same
Painlev\'e--Fredholm hierarchies in a more algebraic form.

\subsection{Resurgent and spectral extensions}

Beyond the polynomial regime, one can pursue a resurgent extension in
which $\zeta$--regularized determinants and Painlev\'e~$\tau$--functions
are regarded as resurgent trans--series\cite{Galapon2022,Rivoal2024}.
The Bernoulli coefficients then appear as Stokes data governing the
connection problem between different asymptotic sectors.
This viewpoint naturally relates the present framework to
the Écalle resurgence program and to the exact~WKB analysis of
Schrödinger operators with analytic potentials.

\begin{figure}[ht!]
\centering
\begin{tikzpicture}[>=Stealth,thick,scale=1.05]
\node (B) at (0,0) {Bernoulli / Stirling};
\node (F) at (4,0) {Fredholm / Painlev\'e};
\node (H) at (2,-1.2) {Hermite / Umbral};
\node (R) at (2,1.2) {Resurgent / WKB};
\draw[->] (B) -- (F);
\draw[->] (B) -- (H);
\draw[->] (H) -- (F);
\draw[dashed,->] (B) -- (R);
\draw[dashed,->] (R) -- (F);
\end{tikzpicture}
\caption{Possible extensions of the Fredholm--Bernoulli bridge:
algebraic (umbral), spectral (Hermite), and asymptotic (resurgent)
directions.}
\end{figure}
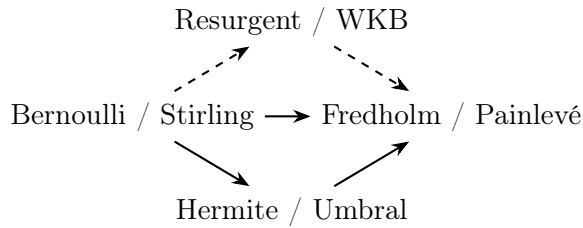

\subsection{Numerical and computational perspectives}

On the computational side, the selector kernels
$\mathcal S_J(\theta,\varphi)$ provide a convenient laboratory for
testing $\zeta$--regularized spectral sums.
Their determinants can be evaluated numerically with high precision
and compared with the asymptotic series
involving $\zeta(2m)$ and $\zeta'(-1)$.
Such comparisons may yield efficient regularization techniques for
finite--rank approximants of integrable kernels and could be extended
to other classes (Airy, Bessel, or sinh kernels).

\subsection{Closing remarks}

The Fredholm--Bernoulli correspondence offers a unified analytic
language for discrete combinatorics, spectral determinants, and
integrable hierarchies.
Its simplicity lies in the observation that
\emph{every regularization is an analytic continuation}.
Through this lens, the Bernoulli numbers serve as
infinitesimal regulators, the $\zeta$--function provides their global
extension, and the Fredholm determinant acts as the functional envelope
that binds them into a single analytic geometry.

\begin{center}
\textit{
From the arithmetic of $B_n$
to the analysis of $\zeta(s)$
to the geometry of Fredholm determinants,
the same continuation principle governs all.}
\end{center}

\vspace{1em}
\noindent
This note has aimed to outline the skeleton of that principle.
Its details---umbral translations, analytic Bernoulli functions,
and resurgent completions---will be developed in subsequent works of
the \emph{hoge \& fuga} series.

\section*{Acknowledgments}
The authors thank the collaborative discussions and numerical verifications
carried out within the \emph{hoge \& fuga} series (2025).


\end{document}